\documentclass[]{article}

\usepackage{lipsum}
\usepackage[a4paper,margin=3.5cm]{geometry}
\usepackage{amsfonts}
\usepackage{graphicx}
\usepackage{epstopdf}
\usepackage{amsmath}          
\usepackage{xfrac}            
\usepackage{bm}               
\usepackage{booktabs}
\usepackage{color}
\usepackage{xcolor}
\usepackage{caption}
\usepackage{subcaption}
\usepackage[
  backend=biber,
  style=numeric,
  giveninits=true]{biblatex}

\bibliography{./bibliography}

\providecommand{\keywords}[1]
{
  \def \and {$\; \cdot \;$}
  \\
  \small
  \noindent
  \textbf{Keywords } {#1}
}

\newenvironment{acknowledgements}
{
 \noindent
 \textbf{Acknowledgments }
}
{
  
}
\begin{document}
\title{Data-driven Stabilization of Nitsche's Method \thanks{Financial support was provided by the German Research Foundation ({\it Deutsche Forschungsgemeinschaft, DFG}) in the framework of the collaborative research center SFB 837 {\it Interaction Modeling in Mechanized Tunneling}.}
}

\author{M. Saberi \thanks{High Performance Computing, Ruhr University Bochum, Universit\"{a}tsstr. 150, 44801 Bochum, Germany} \and
  L. Zhao \footnotemark[2] \and
  A. Vogel \footnotemark[2]}

\maketitle

\begin{abstract}
The weak imposition of essential boundary conditions is an integral aspect of unfitted finite element methods, where the physical boundary does not in general coincide with the computational domain. In this regard, the symmetric Nitsche's method is a powerful technique that preserves the symmetry and variational consistency of the unmodified weak formulation. The stabilization parameter in Nitsche's method plays a crucial role in the stability of the resultant formulation, whose estimation is computationally intensive and dependent on the particular cut configuration using the conventional eigenvalue-based approach. In this work, we employ as model problem the finite cell method in which the need for the generation of a boundary-conforming mesh is circumvented by embedding the physical domain in a, typically regular, background mesh. We propose a data-driven estimate based on machine learning methods for the estimation of the stabilization parameter in Nitsche's method that offers an efficient constant-complexity alternative to the eigenvalue-based approach independent of the cut configuration. It is shown, using numerical benchmarks, that the proposed method can estimate the stabilization parameter accurately and is by far more computationally efficient. The data-driven estimate can be integrated into existing numerical codes with minimal modifications and thanks to the wide adoption of accelerators such as GPUs by machine learning frameworks, can be used with virtually no extra implementation cost on GPU devices, further increasing the potential for computational gains over the conventional eigenvalue-based estimate. The proposed model is tested on both Intel CPU as well as NVIDIA GPU hardware, showing that while it is already many times more efficient on the CPU compared to the eigenvalue-based estimate, its efficiency margin is even larger on modern GPU devices.
\keywords{Nitsche's method \and Finite cell method \and Machine learning \and Data-driven stabilization \and unfitted finite element method}
\end{abstract}
\section{Introduction}
\label{sec:introduction}
The unfitted finite element methods are a class of numerical methods for the approximation of partial differential equations (PDE), whose main aim is to circumvent the need for the generation of a boundary-conforming mesh in conventional finite element methods. The extended finite element method (XFEM)~\cite{moes_1999,belytschko_2001}, the cut finite element method (CutFEM)~\cite{burman_2015}, immersed boundary methods~\cite{peskin_2002,mittal_2005} and the finite cell method (FCM)~\cite{parvizian_2007,duster_2008} fall under the umbrella of unfitted finite element methods. In this work, we focus on the finite cell method, in which the physical domain is extended by a fictitious part and subsequently recovered by means of penalization and a sufficiently accurate integration technique. The extended domain then acts as the computational domain, for whose discretization a regular mesh can be employed, which due to its regular nature allows for efficient computation. The finite cell method has been used in a wide variety of applications, including elasto-plasticity~\cite{abedian_2013_finite,abedian_2014,ranjbar_2014}, linear and nonlinear fluid flow~\cite{saberi_2023_adaptive}, thermo-elasticity~\cite{zander_2012}, shell analysis~\cite{rank_2012}, wave propagation~\cite{duczek_2014,joulaian_2014}, geometrically nonlinear analysis~\cite{schillinger_2012}, topology optimization~\cite{parvizian_2012_topology}, etc., see also~\cite{schillinger_2015} for a review.

As the boundary of the physical domain does not, in general, match the background computational mesh in an embedded setting, boundary conditions are often imposed weakly by means of special techniques. To this end, Lagrange multipliers~\cite{fernandez_2004,glowinski_2007,flemisch_2007,burman_2010}, penalty methods~\cite{babuvska_1973,zhu_1998}, the ghost penalty method~\cite{burman_2010_ghost} and Nitsche's method~\cite{nitsche_1971,hansbo_2002,burman_2012,dolbow_2009,embar_2010} have been used in the literature. In this work, we use the symmetric Nitsche's method. Although penalty-free non-symmetric variants of Nitsche's method have been proposed~\cite{burman_2012_penalty,boiveau_2016_penalty,guo_2017,schillinger_2016,blank_2018}, the major advantage of the original symmetric formulation is the preservation of the symmetry and variational consistency of the unmodified weak formulation. While the weak imposition of boundary conditions in certain techniques, such as the penalty methods, often leads to the numerical ill-conditioning of the resultant system of equations, Nitsche's method is stable and does not suffer from such conditioning issues, given that a sufficiently large stabilization parameter is used, making the choice of the stabilization parameter an important aspect of Nitsche's method, especially for iterative solvers~\cite{saberi_2023_influence}. The stabilization parameter can be estimated using generalized eigenvalue problems, formed either locally for every cutcell or globally for the entire domain~\cite{griebel_2003,dolbow_2009,embar_2010,garhuom_2022,saberi_2023_influence}. The local approach is more computationally efficient, allows for the concurrent computation of the stabilization parameter on each cutcell and was shown to be more robust for multi-level iterative solvers~\cite{embar_2010,saberi_2023_influence}. Nevertheless, the formation and solution of the generalized eigenvalue problems is a computationally costly task, which, as explained in the following sections, requires the integration of multiple terms over the volume of the cutcell as well as across the physical interface with sufficient accuracy. Additionally, given that the resultant matrices are often rank deficient, special techniques for the solution of the generalized eigenvalue problem are necessary, further increasing the associated computational cost.

In this work, we focus on the finite cell method with adaptive mesh refinement (AMR) and adaptive integration for the recovery of the physical domain and Nitsche's method for the imposition of boundary conditions. We use the Poisson equation as model problem. We propose a data-driven approach for the estimation of the stabilization parameter in Nitsche's method and highlight the efficiency of the proposed approach from a computational perspective. The main contributions of this work can be summarized as follows:
\begin{itemize}
\item A data-driven approach based on machine learning methods for the estimation of the stabilization parameter in Nitsche's method is proposed
\item A representation of the cut configuration in the finite cell method with appropriate normalization techniques is presented, and various aspects of designing and training the underlying neural network for the given estimation problem is discussed
\item The proposed data-driven estimate is shown to be on par with the classical eigenvalue estimate in terms of accuracy using a relatively small network
\item The computational efficiency of the proposed method is highlighted by means of numerical benchmarks. It is shown that while the computational cost associated with the estimation of the stabilization parameter using the conventional eigenvalue approach increases as the cut fraction becomes smaller and is, therefore, dependent on the cut configuration, the data-driven approach offers a far more efficient alternative with $O(1)$ computational complexity
\end{itemize}

The remainder of this work is organized as follows. The finite cell formulation of the model problem is described in Section~\ref{sec:fcm}, and the stabilization of Nitsche's method is discussed in Section~\ref{sec:stabilization}. The data-driven stabilization method is presented in Section~\ref{sec:data_driven}. The computational aspects of the proposed approach in comparison with the conventional eigenvalue-based estimate are discussed in Section~\ref{sec:numerical_experiments} using a number of numerical experiments. Finally, conclusions are drawn in Section~\ref{sec:conclusions}.

\section{Model problem}
\label{sec:model_problem}
\subsection{Finite cell formulation}
\label{sec:fcm}
In this work, we use the Poisson equation as model problem and start by deriving its finite cell formulation. The strong form of the Poisson equation can be written as
\begin{equation}
\label{eq:strong}
\begin{split}
- \nabla^{2} u = f \quad\quad &\text{in} \; \Omega,\\
u = g \quad\quad &\text{on} \; \Gamma_{D},\\
\nabla u \cdot \bm{n} = h \quad\quad &\text{on} \; \Gamma_{N}, \\
\end{split}
\end{equation}
where $u$ is the scalar solution variable, $\Omega$ is the spatial domain whose boundary is denoted by $\Gamma = \Gamma_{D} \cup \Gamma_{N}$. $\Gamma_{D}$ and $\Gamma_{N}$ are the Dirichlet and Neumann parts of the boundary, respectively, and $\bm{n}$ is the unit-length outer normal vector to the boundary. $g$ and $h$ are prescribed functions on the Dirichlet and Neumann parts of the boundary, respectively. $f$ is the source term.

We derive the finite cell formulation of the model problem, whose first step, namely multiplication of the strong form with appropriate test functions, integration over the domain $\Omega$ and transferring the derivatives using Green's theorem, is shared with conventional finite element methods, see, e.g.,~\cite{ern_2004}. The boundary-conforming weak form is obtained as:

Find $u \in V$ such that for all $v \in V_{0}$
\begin{equation}
\label{eq:weak_conforming}
\begin{split}
\int_{\Omega} \nabla v \cdot \nabla u \; d \bm{x} - \int_{\Gamma_{D}} v (\nabla u \cdot \bm{n}) \; d \bm{s} = \int_{\Omega} v f \; d \bm{x} +\int_{\Gamma_{N}} v h \; d \bm{s},\\
\end{split}
\end{equation}
where $v$ are the test functions and
\begin{equation}
\begin{split}
V &:= \{u \in H^{1}(\Omega) \; | \; u = g \; \text{on} \; \Gamma_{D}\},\\
V_{0} &:= \{v \in H^{1}(\Omega) \; | \; v = 0 \; \text{on} \; \Gamma_{D}\},\\
\end{split}
\end{equation}
where $H^{1}$ is the Sobolov space.

The domain $\Omega$ is extended to $\Omega_{e}$ by the fictitious part $\Omega_{e} \setminus \Omega$ in the finite cell formulation of the model problem, see Figure~\ref{fig:benchmark}. Therefore, the weak form is modified by a penalization factor $\alpha$, as shown in Equation~\eqref{eq:weak_fcm}, with the aim of recovering the physical domain $\Omega$. Furthermore, whereas the boundary term in Equation~\eqref{eq:weak_conforming} vanishes in the boundary-conforming case as the essential boundary conditions are included in the Space $V$, the finite cell weak formulation with Nitsche's method includes such terms across the physical boundary $\Gamma$, which is not guaranteed to conform to the computational domain $\Omega_{e}$. The weak formulation, after the application of the penalization factor and the addition of Nitsche's method, takes the following form:

Find $u \in V_{e}$ such that for all $v \in V_{e}$
\begin{equation}
\label{eq:weak_fcm}
\begin{split}
\int_{\Omega_{e}} \alpha \nabla v \cdot \nabla u \; d \bm{x} - \int_{\Gamma_{D}} v (\nabla u \cdot \bm{n}) \; d \bm{s} - \int_{\Gamma_{D}} (u - g) (\nabla v \cdot \bm{n}) \; d \bm{s}\\
+ \int_{\Gamma_{D}} \lambda v (u - g) \; d \bm{s} = \int_{\Omega_{e}} \alpha v f \; d \bm{x} +\int_{\Gamma_{N}} v h \; d \bm{s},\\
\end{split}
\end{equation}
where $V_{e} := \{v \in H^{1}(\Omega_{e}) \}$, $\lambda$ is the scalar stabilization parameter in Nitsche's method and the penalization factor $\alpha$ is defined as
\begin{equation}
\label{eq:alpha}
\begin{cases}
\alpha = 1, \quad\quad &\text{ in}\;\Omega,\\
\alpha = 0, \quad\quad &\text{ in}\;\Omega_e \setminus \Omega.
\end{cases}
\end{equation}
In practice, the penalization factor is chosen as $0 < \alpha \ll 1$ instead of zero in the fictitious part $\Omega_{e} \setminus \Omega$ in order to avoid severe numerical ill-conditioning issues. The third and fourth terms on the left-hand side of Equation~\eqref{eq:weak_fcm} are the symmetric consistency and stabilization terms, respectively. The symmetric consistency term ensures that the symmetry of the original formulation is preserved. The imposition of the boundary conditions as well as the stability of the weak form are ensured through the stabilization term, given that a sufficiently large stabilization parameter $\lambda$ is utilized.

Let $\mathcal{T}_{h} := \{K_i\}_{i=1}^{n_K}$ be a tessellation of $\Omega_{e}$ into $n_{K}$ compact, connected, Lipschitz sets $K_{i}$ with non-empty interior and $\mathring{K}_i \cap \mathring{K}_j = \emptyset \; \forall \; i \neq j$. An approximation of the computational domain $\Omega_{e}$ is then defined by $\overline{\Omega}_{e,h} := \cup_{i=1}^{n_{K}} K_{i}$. $\mathcal{T}_{h}$ is typically chosen to be regular in finite cell applications, such as its generation and manipulation can be performed efficiently. We use space tree data structures in this work, see, e.g.,~\cite{burstedde_2011,saberi_2020}, for the discretization of $\Omega_{e}$. We use a combination of adaptive mesh refinement towards $\Gamma$ as illustrated in Figure~\ref{fig:benchmark} and adaptive quadrature integration for the accurate recovery of the physical domain, see~\cite{saberi_2020}. Introducing a finite-dimensional function space $V_{e, h} \subset H^{1}(\Omega_{e})$, the discrete weak form is obtained as follows:

Find $u_{h} \in V_{e, h}$ such that for all $v_{h} \in V_{e, h}$
\begin{equation}
\label{eq:weak_disc_1}
a_{h}(u_{h},v_{h}) = b_{h}(v_{h}),
\end{equation}
with
\begin{equation}
\begin{split}
\label{eq:weak_disc_2}
a_{h}(u_{h},v_{h}) := &\int_{\Omega_{e,h}} \alpha \nabla v_{h} \cdot \nabla u_{h} \; d \bm{x} - \int_{\Gamma_{D,h}} v_{h} (\nabla u_{h} \cdot \bm{n}) \; d \bm{s}\\
&- \int_{\Gamma_{D,h}} u_{h} (\nabla v_{h} \cdot \bm{n}) \; d \bm{s} + \int_{\Gamma_{D,h}} \lambda v_{h} u_{h} \; d \bm{s}, \\[12pt]
b_{h}(v_{h}) := &\int_{\Omega_{e,h}} \alpha v_{h} f_{h} \; d \bm{x} +\int_{\Gamma_{N,h}} v_{h} h_{h} \; d \bm{s} - \int_{\Gamma_{D,h}} g_{h} (\nabla v_{h} \cdot \bm{n}) \; d \bm{s}\\
&+ \int_{\Gamma_{D,h}} \lambda v_{h} g \; d \bm{s},
\end{split}
\end{equation}
where $\Gamma_{D,h}$ and $\Gamma_{N,h}$ are appropriate discretizations of the Dirichlet and Neumann parts of the boundary, respectively.
\subsection{Stabilization}
\label{sec:stabilization}
In this section, we focus on the stabilization of the finite cell weak formulation with Nitsche's method. The symmetric Nitsche's method is stable and does not suffer from numerical ill-conditioning, provided that the stabilization parameter $\lambda$ is sufficiently large. Nevertheless, excessive overestimation of the stabilization parameter leads to numerical ill-conditioning of the resultant system of equations. Therefore, an accurate estimation of the stabilization parameter is an important aspect of Nitsche's method. The estimation of the stabilization parameter in Nitsche's method through the solution of general eigenvalue problems has been studied in the literature, see, e.g.,~\cite{griebel_2003,dolbow_2009,garhuom_2022,saberi_2023_influence}. The bilinear form in Equation~\eqref{eq:weak_disc_2} can be estimated using Young's inequality with $\varepsilon$ as:
\begin{equation}
\label{eq:weak_stabilization_1}
\begin{split}
a_{h}(v_{h},v_{h}) \geq &\int_{\Omega_{e,h}} \alpha \nabla v_{h} \cdot \nabla v_{h} \; d \bm{x} - \frac{1}{\varepsilon} \int_{\Gamma_{D,h}} v_{h} v_{h} \; d \bm{s}\\
&- \varepsilon \int_{\Gamma_{D,h}} (\nabla v_{h} \cdot \bm{n})(\nabla v_{h} \cdot \bm{n}) \; d \bm{s} + \int_{\Gamma_{D,h}} \lambda v_{h} v_{h} \; d \bm{s}, \; \forall \varepsilon > 0. \\
\end{split}
\end{equation}
Assuming the stabilization parameter $\lambda$ is constant over the integration domain and assuming a constant scalar $C$ such that
\begin{equation}
\label{eq:weak_stabilization_2}
\begin{split}
C \int_{\Omega_{e,h}} \alpha \nabla v_{h} \cdot \nabla v_{h} \; d \bm{x} \geq \int_{\Gamma_{D,h}} (\nabla v_{h} \cdot \bm{n})(\nabla v_{h} \cdot \bm{n}) \; d \bm{s},\\
\end{split}
\end{equation}
the following inequality can be obtained from Inequality~\eqref{eq:weak_stabilization_1}:
\begin{equation}
\label{eq:weak_stabilization_3}
\begin{split}
a_{h}(v_{h},v_{h}) \geq (1 - \varepsilon C)  \int_{\Omega_{e,h}} \alpha \nabla v_{h} \cdot \nabla v_{h} \; d \bm{x} + (\lambda - \frac{1}{\varepsilon}) \int_{\Gamma_{D,h}} v_{h} v_{h} \; d \bm{s}.\\
\end{split}
\end{equation}
$(1 - \varepsilon C)$ and $(\lambda - \frac{1}{\varepsilon})$ must be both positive for the coercivity of the bilinear form. It follows that $\lambda > C$. Thereby, a lower bound for $\lambda$ is obtained. The constant $C$ can be computed as the largest eigenvalue of the following generalized eigenvalue problem:
\begin{equation}
\label{eq:gen_eig_1}
\mathbf{K} \bm{v} = \Lambda \mathbf{M} \bm{v},
\end{equation}
where $\Lambda$ and $\bm{v}$ are the eigenvalues and eigenvectors, respectively. The matrices $\mathbf{K}$ and $\mathbf{M}$ are given by
\begin{equation}
\label{eq:gen_eig_2}
\begin{cases}
\mathbf{K}_{ij} :=  \int_{\Gamma_{D,h}} (\nabla \phi_{i} \cdot \bm{n})(\nabla \phi_{j} \cdot \bm{n}) \; d \bm{s},\\[12pt]
\mathbf{M}_{ij} := \int_{\Omega_{e}} \alpha \nabla \phi_{i} \cdot \nabla \phi_{j} \; d\bm{x},
\end{cases}
\end{equation}
and using a basis $\{ \phi_i \} \subset V_{e,h}$. Therefore, the lower bound for the stabilization parameter in Nitsche's method can now be chosen as the largest eigenvalue $\max_{k}{\Lambda_{k}}$. Although the terms in Inequalities~\eqref{eq:weak_stabilization_1} -~\eqref{eq:weak_stabilization_3} and, consequently, matrices $\mathbf{M}$ and $\mathbf{K}$ are integrated over the entire computational domain $\Omega_{e,h}$ and the Dirichlet boundary $\Gamma_{D,h}$, respectively, the integration domain can be restricted to the domain of a single cell $K$ conservatively, thereby effectively leading to a local estimate that can be computed concurrently for every $K \in \mathcal{T}_{h}$. The matrices $\mathbf{M}$ and $\mathbf{K}$ in Equation~\eqref{eq:gen_eig_1} can then be computed as follows for every cell $K$:
\begin{equation}
\label{eq:gen_eig_3}
\begin{cases}
\mathbf{K}_{ij} :=  \int_{\Gamma_{D,h}^{K}} (\nabla \phi_{i} \cdot \bm{n})(\nabla \phi_{j} \cdot \bm{n}) \; d \bm{s},\\[12pt]
\mathbf{M}_{ij} := \int_{K} \alpha \nabla \phi_{i} \cdot \nabla \phi_{j} \; d\bm{x},
\end{cases}
\end{equation}
where $\Gamma_{D,h}^{K}$ is the Dirichlet boundary of cell $K$. The stabilization parameter can now be chosen as a cell-wise constant, computed as the largest eigenvalue of the local generalized eigenvalue problem in Equation~\eqref{eq:gen_eig_3}.

The solution to the generalized eigenvalue problem in Equation~\eqref{eq:gen_eig_1} is in general non-trivial due to the rank deficiency of the matrices, and typically scales poorly with problem size; therefore, aside from a parallel computing standpoint, the local estimate, because of its smaller size, is more favorable also from the perspective of computational complexity. Furthermore, it was shown in~\cite{saberi_2023_influence} that such local estimation is advantageous in the context of multi-level solvers, specifically the geometric multigrid method. Therefore, we focus on the local estimate in this work. We note that the stabilization parameter obtained from the solution of the eigenvalue problem above is the lower bound for the stability of the formulation. In practice, an additional safety factor, typically chosen as 2, is used for the computed stabilization parameter.

\section{Data-driven stabilization}
\label{sec:data_driven}
In this section, we propose a data-driven approach for the approximation of the stabilization parameter in Nitsche's method. The stabilization parameter $\lambda$ obtained from the generalized eigenvalue problem in Equations~\eqref{eq:gen_eig_1} and~\eqref{eq:gen_eig_3} can be represented as a function $\lambda = G(\mathbf{x})$, in which $\mathbf{x}$ is a vector representing the local cut configuration in a given cell. We propose to approximate the function $G(\mathbf{x})$ with the help of a neural network. To this end, the cut configuration, i.e., the local boundary $\Gamma_{K}$ within each cell, must be represented as a real-valued vector $\mathbf{x}$ in order to be used as the input to the neural network, whose selection is explained in detail in Section~\ref{sec:cut_config}. In the following, we outline the design and performance of the proposed data-driven approach to the approximation of the stabilization parameter. We focus on the 2D case in this work to establish the proposed concept; nevertheless, we note that similar steps can in principle be applied to the 3D case.

\subsection{Representation of the cut configuration}
\label{sec:cut_config}
\begin{figure}[tb]
    \centering
    \includegraphics{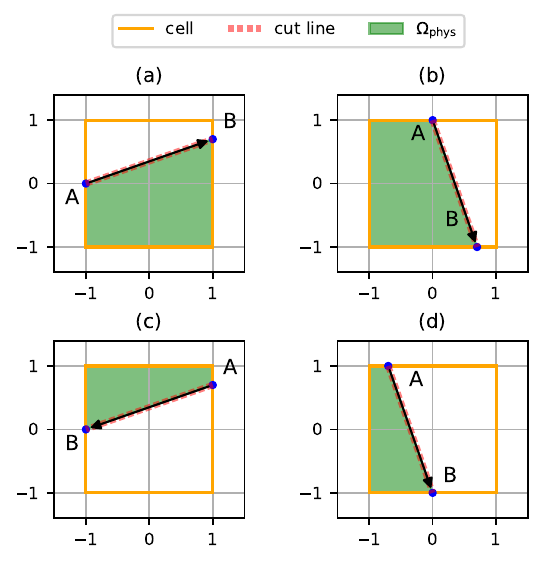}
    \caption{The representation of the cut configuration and its normalization. (a), (c) two cut configurations before normalization with the same cut line and flipped starting and end points. (b), (d) the normalized form of the cut configurations in (a) and (c), respectively. Note that the normalized representation can distinguish between the cut configurations in (a) and (c)}
    \label{fig:cut_config}
\end{figure}
The approximation of $G(\mathbf{x})$ using a neural network presupposes the existence of an appropriate representation of the local physical boundary $\Gamma_{K}$ in the form of a real-valued vector $\mathbf{x}$ that can describe the essential features of $\Gamma_{K}$ sufficiently accurately. As the stabilization parameter, computed from the eigenvalue problem in Equations~\eqref{eq:gen_eig_1} and~\eqref{eq:gen_eig_3}, depends on the location of the local boundary, it is incumbent on such representation to provide a sufficiently sensitive mapping between the cut configuration and the corresponding stabilization parameter, and consequently, the performance of the network heavily depends on the choice of this representation. We propose a representation technique for the local boundary in the following and discuss its selection.

The cut configuration is approximated using a line $L$, referred to as the cut line and defined by the intersection points $A$ and $B$ between the physical boundary and the boundaries of the cell, the justification for whose assumption is twofold. First, practical observations indicate that even when the physical boundary $\Gamma$ is curved, the local boundary $\Gamma_{K}$ is, in the vast majority of cases, a line or closely approximated by a line when the computational domain $\Omega$ is adaptively refined towards the physical boundary, as described in Section~\ref{sec:fcm}. This is in particular true if the boundary is described in terms of manifolds, e.g., polygons, splines, etc. This point is further illustrated in Section~\ref{sec:numerical_experiments}. Second, deviations of $\Gamma_{K}$ from a line are easy to detect and quantify, allowing the computation of the stabilization parameter to fall back to the generalized eigenvalue problem in exceptional cases. For instance, the ratio between the length of $L$ and the length of $\Gamma_{K}$ can act as a simple metric to quantify how well the latter is approximated by the former. Furthermore, although the physical boundary may cut a given cell into multiple regions, such cases are even rarer in practice.

It is important to note that the physical boundary divides a cutcell into a physical and a fictitious region, which must be distinguished by the representation of the cut configuration; therefore, the cut line is treated as a vector, defined by the points $A$ and $B$, as shown in Figure~\ref{fig:cut_config}, where the physical domain $\Omega_{\text{phys}}$ is on the right side of $\vec{AB}$.

We propose to represent the cut configuration as the distance between the line $L$ and $n$ predetermined feature points $T_{i}, \; i = 1, \cdots, n$, referred to as the cut distance. We arrange the cut distances into a vector $\mathbf{x}$ and use it as the input to the neural network such that $x_{i}$ corresponds to the distance between $L$ and $T_{i}$. The computational cost associated with the preparation of the vector $\mathbf{x}$ is low as the calculation of the distance between a line and a point is a relatively cheap operation. Furthermore, such representation allows the neural network to operate within the logarithmic space as the distances are always non-negative, the importance of which is discussed in Section~\ref{sec:dataset}. More importantly, the proposed representation allows for the flexible adjustment of the input features through the feature points. It should be noted that the vector $\mathbf{x}$ does not contain any information regarding the direction of the cut line $L$; therefore, we propose a normalization technique for the implicit inclusion of such information. More specifically, we select the top edge of the cell as the starting edge, on which the starting point of the cut vector $\vec{AB}$ lies. Since the generalized eigenvalue problem in Equation~\eqref{eq:gen_eig_1}, and therefore, the stabilization parameter $\lambda$ are invariant to the rotation of the standard cell, a simple rotation would be sufficient for such normalization in cases where the start point of $\vec{AB}$ lies on another edge of the cell. The normalization process is demonstrated in Figure~\ref{fig:cut_config}, where two cut configurations with the same cut line but flipped starting and end points in Figure~\ref{fig:cut_config}~(a) and (c) are transformed to their normalized configuration in Figure~\ref{fig:cut_config}~(b) and (d), respectively.

The standard cell, as shown in Figure~\ref{fig:cut_config}, is used to compute the cut distances in order to keep the input space bounded. The feature points should be chosen such that 1) $\mathbf{x}$ contains sufficient information for the accurate characterization of the relationship between the cut configuration and the stabilization parameter by the network and 2) $\mathbf{x}$ is as small as possible in order to limit the associated computational cost. The selection of the feature points is carried out with the help of a hyperparameter search which is explained in Section~\ref{sec:network}.
\subsection{Dataset}
\label{sec:dataset}
\begin{figure}[tb]
    \centering
    \includegraphics{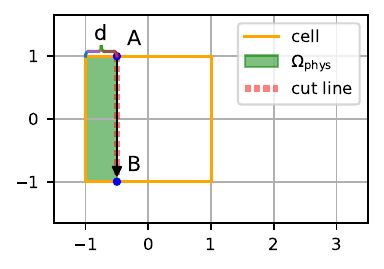}
    \caption{The rectangular cut configuration, where the cut sliver is defined using the width of the rectangle $d$}
    \label{fig:rectangle_cut_config}
\end{figure}
\begin{figure}[tb]
    \centering
    \includegraphics{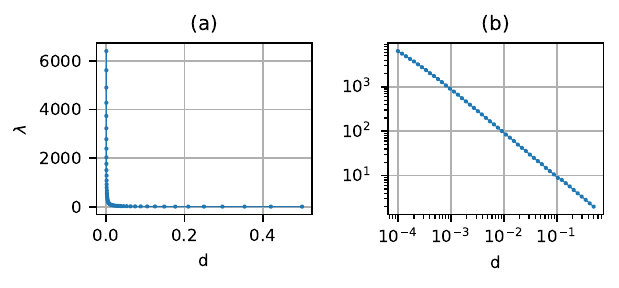}
    \caption{The stabilization parameter from the rectangular cut configuration test, where the width of the rectangular cut sliver $d$, see Figure~\ref{fig:rectangle_cut_config}, is progressively reduced in (a) linear scale and (b) logarithmic scale. The data is presented in order of the values of the stabilization parameter $\lambda$}
    \label{fig:rectangle_cut_res}
\end{figure}
\begin{figure}[tb]
    \centering
    \includegraphics{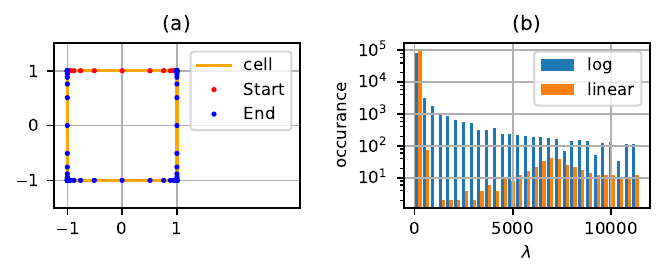}
    \caption{(a) Logarithmic spacing of feature points on the edges of the cell, where the start points are placed on the top edge and the end points are placed on the other edges according to the normalization technique in Section~\ref{sec:cut_config} and (b) the distribution of the stabilization parameter $\lambda$ with linear and logarithmic spacing of start and end points}
    \label{fig:target_points}
\end{figure}
\begin{figure}[tb]
    \centering
    \includegraphics{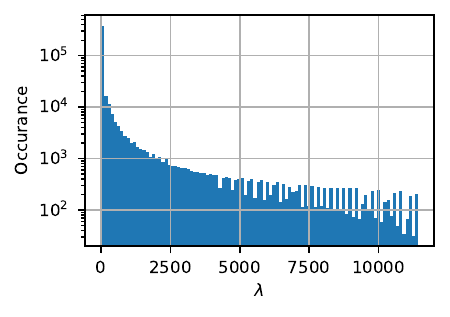}
    \caption{The distribution of the stabilization parameter $\lambda$ in the training dataset, where 399 feature points with logarithmic spacing are placed on each edge of the cell, resulting in a total of $477,603$ data points}
    \label{fig:training_dataset}
\end{figure}
The training dataset must provide a sufficient representation of possible cut configurations. We employ a distribution of points on each edge of the standard cell to generate the dataset by treating them as the start and end points of the cut configuration. In order to better understand the behavior of the objective function $G(\mathbf{x})$, we first compute the stabilization parameter $\lambda$ using the eigenvalue estimate for a vertical cut line, as shown in Figure~\ref{fig:rectangle_cut_config}, where the cut configuration can be defined using a single parameter, namely the width of the cut sliver, $d$. It can be seen in Figures~\ref{fig:rectangle_cut_res}~(a) and (b) that the stabilization parameter $\lambda$ closely follows exponential growth as the cut sliver becomes smaller. In light of this behavior and in order to strike a balance in the distribution of the stabilization parameter in the dataset between benign and severe cut configurations, it can be inferred that the presence of cases where the cut fraction is small must be heightened. To this end, the start and end points of the cut configurations can be distributed logarithmically towards the vertices of the cell such that for a given number of points on any edge, the density of points is exponentially higher near the vertices, as shown in Figure~\ref{fig:target_points}~(a). It can be seen in Figure~\ref{fig:target_points}~(b) that the logarithmic spacing of the start and end points, as expected, leads to a more uniform distribution of the stabilization parameter compared to its linear counterpart. Therefore, the training dataset is generated using 399 points, spaced logarithmically, on each edge, where the minimum distance of $10^{-4}$ to the vertices is chosen such that the dataset contains a realistic representation of cut configurations that the network may encounter. The distribution of the stabilization parameter in the training dataset is shown in Figure~\ref{fig:training_dataset}.

Given the rotational invariance of the stabilization parameter to the cut configuration, as pointed out in Section~\ref{sec:cut_config}, traversing all combinations of start and end points in such setting would lead to duplicated entries in the dataset, which in addition to the wasted computational effort would hinder the training of the network. The normalization technique in Section~\ref{sec:cut_config} makes it possible to eliminate such duplicate cases by treating the points on the top edge as starting points and the points on the remaining edges as end points, and the distribution of the stabilization parameter, shown, e.g., in Figures~\ref{fig:rectangle_cut_res} and~\ref{fig:training_dataset}, does not contain duplicate cases.

The validation and test datasets with 118,803 and 96,123 data points, respectively, are produced in a similar fashion as the training dataset, where overlaps between the datasets are avoided by using different distributions of the points and, consequently, cut configurations.

As shown in Figure~\ref{fig:rectangle_cut_res}, the stabilization parameter grows exponentially as the cut fraction tends to zero. It is well-known that neural networks do not perform well for such highly nonlinear functions; therefore, we transform both the input and output into logarithmic spaces, where the relationship between the cut fraction and the stabilization parameter is almost linear, as shown Figure~\ref{fig:rectangle_cut_res}~(b). Such transformation brings about the need for another normalization step as cut configurations that pass exactly through a feature point $T_{i}$ would invalidate the corresponding $i$-th entry in vector $\mathbf{x}$ as the cut distance would be zero. Therefore, we use a minimum cut-off value of $10^{-10}$ for the cut distances before transforming the input into logarithmic space.
\subsection{Network}
\label{sec:network}
\begin{figure}[tb]
    \centering    
    \includegraphics{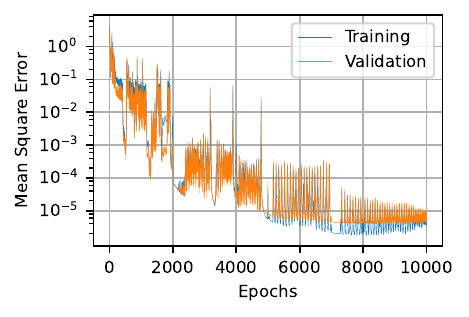}    
    \caption{The training and validation loss of the final network with a depth of 6, a width of 1024 at each layer and the input feature $T_{v} + T|_{0.002}$, see Table~\ref{tab:hyper_search_1}, over 10000 epochs}
    \label{fig:training_loss}
\end{figure}
\begin{figure}[tb]
    \centering
    \includegraphics{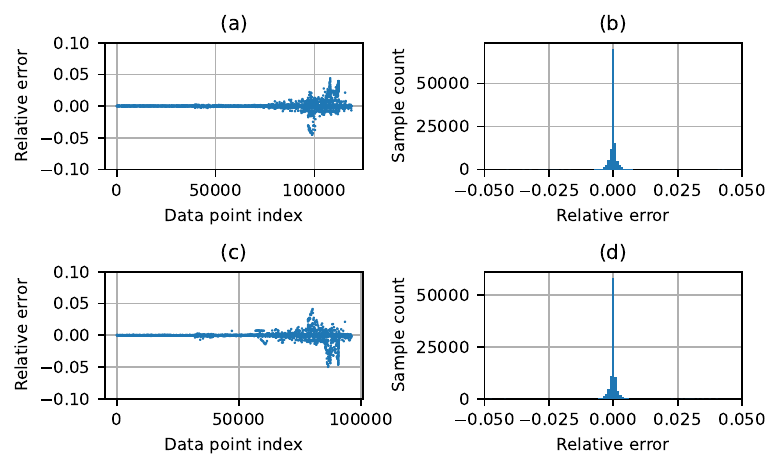}
    \caption{The accuracy of the final network with 6 fully-connected layers, each with 1024 neurons and the input feature $T_{v} + T|_{0.002}$, see Table~\ref{tab:hyper_search_1}. (a), (c) the relative error of the network for the validation and test datasets, respectively, where the data points are sorted by the value of the stabilization parameter. (b), (d) the distribution of the relative error for the validation and test datasets, respectively}
    \label{fig:network_accuracy}
\end{figure}
\begin{table}[tb]
    \centering
    \caption{The configurations of the feature points tested in the first stage of the hyperparameter search. $T_{v}$ denotes the vertices of the cell, $T_{G}^{n}$ denotes the Gauss integration points with degree $n$ in each direction, $T_{\cdot, \cdot}$ denotes a feature point at coordinates $(\cdot, \cdot)$ and $T|_{(\cdot)}$ denotes feature points located with a distance of $(\cdot)$ from each vertex. $n_{T}$ denotes the total number of feature points in a given layout. The layouts with the lowest validation loss are shown in boldface}
    \label{tab:hyper_search_1}
    \begin{tabular}{lll}
    \toprule
    Input feature points & $n_{T}$ & Validation Loss   \\
    \midrule
    $T_{v}$ & 4 & 0.12942 \\
    $T_{v} + T_{(0, 0)}$ & 5 & 0.11291 \\
    $T_{v} + T_{G}^{2}$ & 8 & 0.00592 \\
    $T_{G}^{2}$ & 4 & 0.00583 \\
    $T_{v} + T_{\text{lin}}^{1}$ & 8 & 0.03182 \\
    $T_{v} + T_{\text{lin}}^{2}$ & 12 & 0.01899 \\
    $T_{v} + T_{\text{lin}}^{3}$ & 16 & 0.03338 \\
    $T_{v} + T_{\text{lin}}^{4}$ & 20 & 0.00329 \\
    $T_{v} + T_{\text{lin}}^{5}$ & 24 & 0.00524 \\
    $T_{v} + T_{\text{lin}}^{7}$ & 32 & 0.00533 \\
    $T_{v} + T_{\text{lin}}^{10}$ & 44 & 0.00099 \\
    $T_{v} + T_{\text{lin}}^{15}$ & 64 & 0.00381 \\
    $T_{v} + T_{\text{lin}}^{20}$ & 84 & 0.00088 \\
    $T_{v} + T_{\text{lin}}^{40}$ & 164 & 0.00064 \\
    $T_{v} + T|_{(0.2)}$ & 12 & 0.00058  \\
    $T_{v} + T|_{(0.0632)}$ & 12 & 0.00036  \\
    $T_{v} + T|_{(0.02)}$ & 12 & 0.00037  \\
    $T_{v} + T|_{(0.00632)}$ & 12 & 0.00027  \\
    $T_{v} + T|_{(0.002)}$ & 12 & \textbf{0.00019}  \\
    $T_{v} + T|_{(0.000632)}$ & 12 & 0.00020  \\
    $T_{v} + T|_{(0.0002)}$ & 12 & 0.00097  \\
    $T_{v} + T|_{(0.0000632)}$ & 12 & 0.00454  \\
    $T|_{(0.002)}$ & 8 & 0.00064  \\
    $T_{v} + T|_{(0.2)} + T|_{(0.02)}$ & 20 & \textbf{0.00019} \\
    $T_{v} + T|_{(0.02)} + T|_{(0.002)}$ & 20 & \textbf{0.00019} \\
    $T_{v} + T|_{(0.002)} + T|_{(0.0002)}$ & 20 & 0.00022 \\
    \bottomrule
    \end{tabular}
\end{table}
\begin{table}[tb]
\centering
\caption{The configurations of the hidden layers tested in the second stage of the hyperparameter search. The input feature $T_{v} + T|_{(0.002)}$ from the first stage of the hyperparameter search is used. The two lowest validation losses are shown in boldface}
\label{tab:hyper_search_2}
\begin{tabular}{lll}
    \toprule
    Depth & Width & Validation Loss \\
    \midrule
    5 & 512 & 4.663e-5 \\
    5 & 1024 & 4.249e-5 \\
    5 & 2048 & 3.261e-5 \\
    5 & 4096 & 2.927e-5 \\
    6 & 512 & 3.149e-5 \\
    6 & 1024 & \textbf{2.328e-5} \\
    6 & 2048 & 2.994e-5 \\
    6 & 4096 & 3.230e-5 \\
    7 & 512 & 3.043e-5 \\
    7 & 1024 & 2.237e-5 \\
    7 & 2048 & \textbf{2.318e-5} \\
    7 & 4096 & 5.370e-5 \\
    \bottomrule
\end{tabular}
\end{table}
As the dimension of the input layer, i.e., vector $\mathbf{x}$, is low and the output, i.e., the stabilization parameter $\lambda$, is continuous, we use a fully connected feedforward neural network with ReLu activation function on all layers, see, e.g.,~\cite{goodfellow_2016}. The neural network is implemented using TensorFlow~\cite{tensorflow_2015}. The Adam optimization algorithm is used in all cases. A normalization layer is placed right after the input layer, which applies shifting and scaling to achieve a distribution with zero mean and unit standard deviation. The normalization is applied individually on each element of $\mathbf{x}$, i.e., the dimension of $\mathbf{x}$ is unchanged. As the relative error of the stabilization parameter plays an important role in the stability of the finite cell formulation, we define outliers as predicted output with more than $5\%$ error and use their frequency as a metric for the accuracy of the network. The mean square error (MSE) loss is selected as the loss function, which provides a good indication of the network accuracy as MSE losses are observed to correspond to fewer outliers.

It was found that in addition to the configuration of the hidden layers, the layout of the feature points in the proposed representation of the cut configuration has a decisive impact on the prediction accuracy of the network; therefore, in order to determine the optimum network architecture, i.e., the depth and width of the network as well as the number and location of the feature points at the input layer, a two-stage hyperparameter search is conducted. The mean square error of the network on the validation set is used as the optimization criterion. The hyperparameter search is carried out in two stages mainly to alleviate the associated computational cost. In each trial, the network is trained for 1000 epochs with an exponential decay learning rate schedule. The initial learning rate is $0.0005$ and is halved every 125 epochs. We note that the low initial learning rate and the exponential decay schedule, together with a large batch size of 65,536 contribute to reduce overshoot (oscillating losses) in the training.

During the first stage, the optimum set of feature points is sought between a number of predefined configurations, as shown in Table~\ref{tab:hyper_search_1}, by fixing the hidden layers at $6 \times 512$, which was found to provide reasonable results. As shown in Table~\ref{tab:hyper_search_1}, a wide range of layouts with different number and placement of the feature points are tested. The results indicate that the addition of feature points leads to lower validation loss in general. Furthermore, it can be seen that additional feature points close to the vertices in specific significantly improve the performance of the network, for instance, placing feature points at a distance of $0.002$ from each vertex in addition to the vertices leads to a lower validation loss compared to placing 40 linearly spaced feature points on each edge. Given that, due to the high sensitivity of the stabilization parameter to relative changes to the cut configuration, the most challenging cases for the network to predict are expected to occur where the stabilization parameter is large, i.e., small cut fractions, it is suspected that the observed improvement can be attributed to the fact that the addition of points very close the vertices provides the network with useful information by rendering the input more sensitive to changes precisely in such cases. In addition, the accuracy of the network is relatively robust with respect to the placement of the feature points as long as the additional points are close to the vertices and the validation loss is not significantly affected when the extra feature points are placed between roughly $10^{-1}$ and $10^{-4}$ from the vertices. It was found that the four vertices of the cell along with one point located close to each vertex on each edge result in a good compromise between accuracy and computational cost. We choose the input feature $T_{v} + T|_{0.002}$, see Table~\ref{tab:hyper_search_1}, for the final network. We note that an increase in the number of feature points would entail additional distance computations on each cell, and a lower input dimension is, therefore, highly preferred.

In the second stage of the hyperparameter search, the input layer was fixed using the results from the first stage and the optimum configuration of the hidden layers was sought between the cases shown in Table~\ref{tab:hyper_search_2}. It was found that the network with $6 \times 1024$ neurons in the hidden layers leads to a good compromise between network size and accuracy. Furthermore, it can be seen that the depth and width of the network have a more limited impact on accuracy in comparison with the layout of the feature points. We note that repeating the first stage by fixing the hidden layers at $6 \times 1024$ did not change the outcome of the search in terms of the relative validation loss of the different layouts.

The final network is trained for $10000$ epochs with an exponential decay learning rate schedule and a batch size of 65,536, where the initial learning rate was 0.0005 and halved every 2,500 epochs. The checkpoint feature of TensorFlow is employed in order to preserve the model parameters with the minimum validation loss. The network loss for the training and validation sets is shown in Figure~\ref{fig:training_loss}, where a good convergence can be observed for both datasets. We note that training the network for more epochs did not, in our tests, lead to reliably improved convergence.

The mean square error of the network on the validation set as well as the relative error on individual data points were used to evaluate the performance of the network. As the underestimation of the stabilization parameter may lead to instability of the finite cell discretization, it is important to control the frequency and intensity of outlier data points. The relative error, calculated as $e_{r} := \frac{\lambda_{p} - \lambda_{g}}{\lambda_{g}}$, where $\lambda_{p}$ and $\lambda_{g}$ are the predicted and ground truth stabilization parameters, respectively, are shown in Figure~\ref{fig:network_accuracy}. It can be seen that the maximum relative error is bounded to around $5\%$, which is well within the stable range given that a safety factor of 2 for the stabilization parameter is used in practice, and the vast majority of cases, as seen in Figure~\ref{fig:network_accuracy}~(b) and (d), are well below $2.5\%$.
\subsection{Implementation and integration into existing codes}
A major advantage of the presented data-driven estimate is its straightforward integration into existing simulation codes. The offline phase of the data-driven estimate consists in the training of the underlying neural network using appropriate data and has to be performed once for a given problem type. The online phase of the estimate consists in the integration of the model into the existing simulation pipeline by replacing the routine responsible for the computation of the eigenvalue problem with the data-driven estimate, which can be achieved without intrusive modifications to the code. The input to the data-driven estimate, which is composed of the cut distances on each cutcell, can be obtained either on the fly for each cutcell or in a pre-processing step for all cutcells at once. As shall be seen in Section~\ref{sec:numerical_experiments}, the evaluation of the stabilization parameter using the data-driven estimate, which corresponds to the inference of the neural network, is typically more efficient in batches. Therefore, the latter, i.e., the pre-processing of the input data on all cutcells and batch estimation of the stabilization parameter is preferred in practice. Furthermore, thanks to the provision of interfaces to different programming languages by machine learning frameworks, the integration of the data-driven estimate is typically without barriers. We note that the normalization step in Section~\ref{sec:cut_config} can be circumvented in practice simply by reordering the feature points $T_{i}$ in vector $\mathbf{x}$ as a function of the starting edge and, therefore, does not entail any additional computational effort. As the data-driven estimate is trained using the standard cell, the estimated stabilization parameter must be transformed for a cell of arbitrary size. Specifically, the stabilization parameter scales linearly with the cell size for a square cell; therefore, a scaling factor of $\frac{2}{l}$, where $l$ is the side length of the cell is applied to the output of the data-driven estimate. Finally, as accelerators such as GPUs are supported by most machine learning frameworks, the estimation of the stabilization parameter can be offloaded to such devices without any extra implementation effort using the data-driven estimate.
\section{Numerical experiments}
\label{sec:numerical_experiments}
\begin{figure}[tb]
    \centering
    \includegraphics{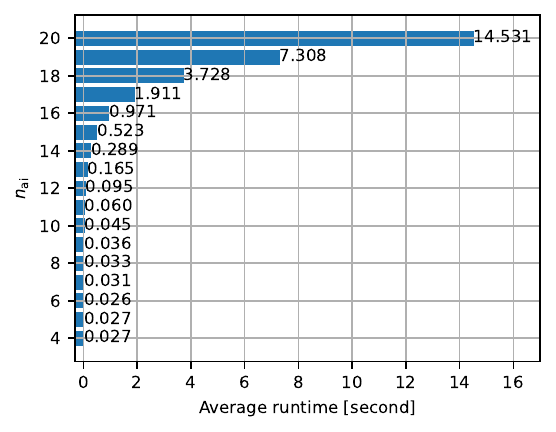}
    \caption{The total runtime of the eigenvalue approach for the computation of the stabilization parameter $\lambda$ on a single cutcell with different number of adaptive integrations $n_{\text{ai}}$ on a single thread of an Intel Xeon E5-2630 v3 CPU running at 2.40 GHz}
    \label{fig:eigval_runtime}
\end{figure}
\begin{figure}[tb]
    \centering
    \includegraphics{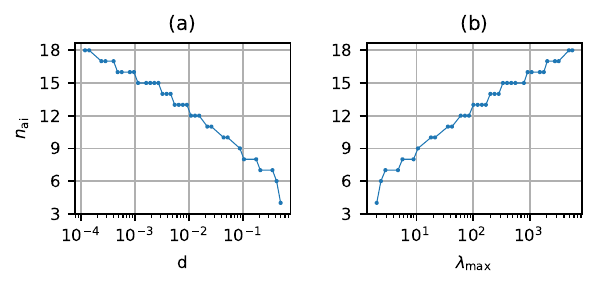}
    \caption{The necessary number of adaptive integration levels $n_{\text{ai}}$ as a function of (a) the cut sliver and (b) the maximum eigenvalue in the vertical cut test in Figure~\ref{fig:rectangle_cut_config}. The required $n_{\text{ai}}$ is determined by calculating $\lambda$ starting with $n_{\text{ai}} = 3$ and progressively increasing $n_{\text{ai}}$ until the relative change in $\lambda$ is below $1\%$ of the reference value calculated with $n_{\text{ai}} = 20$}
    \label{fig:nai}
\end{figure}
\begin{figure}[tb]
    \centering
    \includegraphics{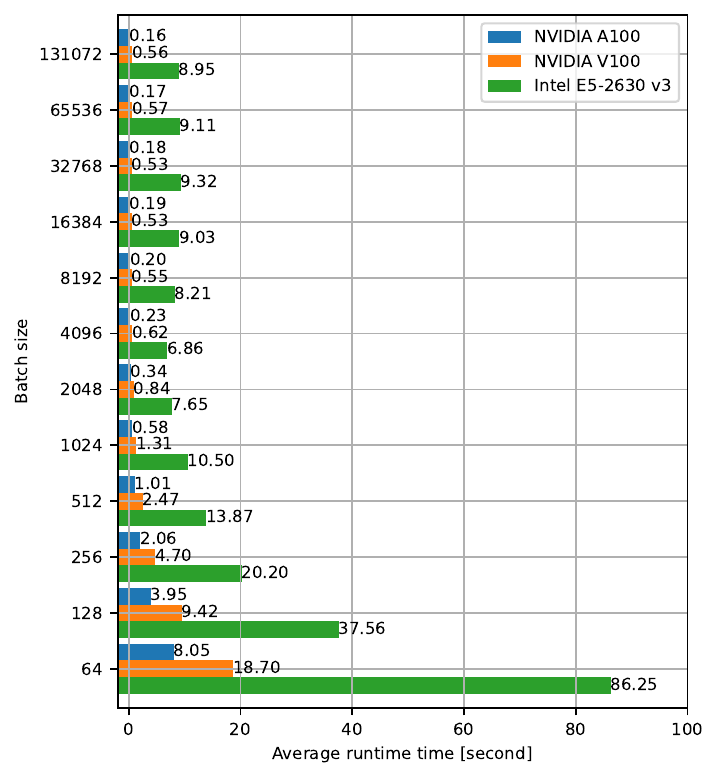}
    \caption{The runtime of the data-driven approach for the estimation of the stabilization parameter $\lambda$ on the entire training dataset with 477,603 data points using the Intel Xeon E5-2630 v3, the NVIDIA Tesla V100 and the NVIDIA Ampere A100 using different batch sizes}
    \label{fig:performance_training}
\end{figure}
\begin{figure}[tb]
    \centering
    \includegraphics{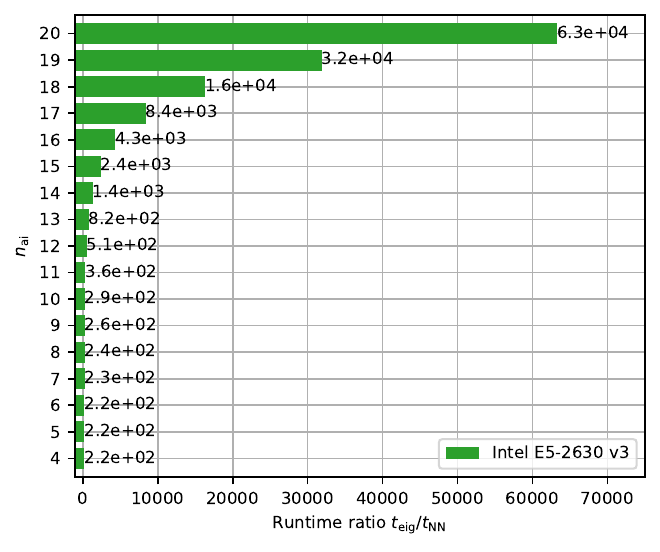}
    \caption{The runtime ratio between the eigenvalue approach and the data-driven approach on the Intel Xeon E5-2630 v3 CPU for the entire training dataset}
    \label{fig:runtime_ratio_cpu}
\end{figure}
\begin{figure}[tb]
    \centering
    \includegraphics{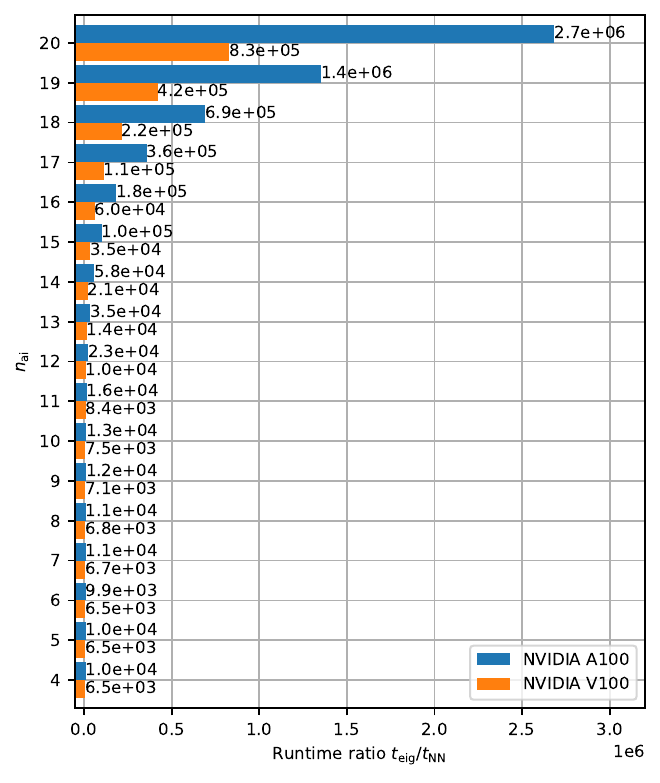}
    \caption{The runtime ratio between the eigenvalue approach on the Intel Xeon E5-2630 v3 CPU and the data-driven approach on the NVIDIA Tesla V100 and NVIDIA Ampere A100 GPUs for the entire training dataset}
    \label{fig:runtime_ratio_gpu}
\end{figure}
\begin{figure}[tb]
    \centering
    \includegraphics{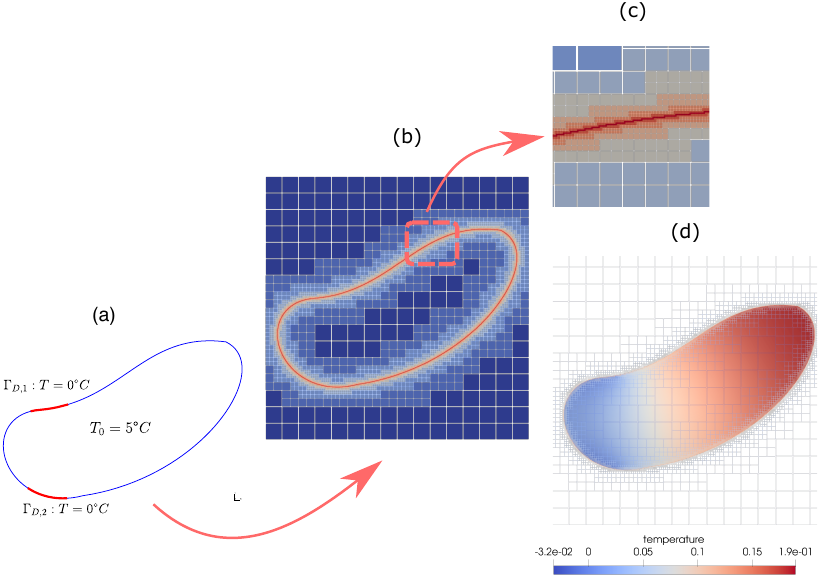}
    \caption{The numerical benchmark for the evaluation of the data-driven approach for the estimation of the stabilization parameter. (a) The physical domain and boundary conditions, (b) the discretization of the computational domain with adaptive mesh refinement towards the boundary, (c) a close-up look at a part of the computation domain around the physical boundary and (d) the temperature field on the physical boundary}
    \label{fig:benchmark}
\end{figure}
In this section, the performance of the proposed data-driven estimate in terms of accuracy and computational efficiency is studied using a number of benchmarks. The numerical experiments are carried out using an in-house C++ implementation, where \texttt{p4est}~\cite{burstedde_2011} and \texttt{PETSc}~\cite{petsc_1997} are used for mesh manipulation and some linear algebra operations, respectively. \texttt{Paraview}~\cite{ahrens_2005} is used for the visualization and postprocessing of some of the results.

The calculation of the stabilization parameter $\lambda$ via the generalized eigenvalue problem in Equation~\eqref{eq:gen_eig_1} is computationally expensive mainly on account of the high cost associated with the computation of the matrices $\mathbf{K}$ and $\mathbf{M}$ on the one hand and the solution of the generalized eigenvalue problem on the other. In this regard, as the total effort for the integration of $\mathbf{K}$ and $\mathbf{M}$ scales with the number of quadrature points, the number of adaptive integration levels $n_{\text{ai}}$ directly affects the associated computational cost. Furthermore, the solution of the resultant generalized eigenvalue problem is similarly expensive, especially given that the matrices $\mathbf{K}$ and $\mathbf{M}$ are typically rank deficient. The average total runtime for the calculation of the stabilization parameter on a single cutcell using a single thread of an Intel Xeon E5-2630 v3 processor is shown in Figure~\ref{fig:eigval_runtime}. It can be seen that, as expected, the runtime significantly grows with the number of adaptive integration levels. Furthermore, depending on the implementation, the memory footprint of such computation, which similarly grows with the number of adaptive integration levels, can become a bottleneck, especially given that the stabilization parameter for multiple cells can be processed in parallel. The situation is exacerbated by the fact that a priori determination of an appropriate $n_{\text{ai}}$ for a given cut configuration is not trivial, and a brute-force approach would either entail the repeated computation of the stabilization parameter with increasing $n_{\text{ai}}$ until a converged solution is obtained or the selection of a large $n_{\text{ai}}$ for all cases, both of which are computationally wasteful. Regardless of the lack of knowledge of the optimal $n_{\text{ai}}$ for a given cut configuration, severe cut configuration, i.e., small cut fractions, require a larger number of adaptive integration levels, as shown in Figure~\ref{fig:nai}, and are, consequently, more computationally costly. In Figure~\ref{fig:nai}~(a), the necessary $n_{\text{ai}}$ for the accurate calculation of the stabilization parameter is shown as a function of the cut sliver in Figure~\ref{fig:rectangle_cut_config}, in which the dependence of the required $n_{\text{ai}}$ on the cut size is clearly visible.

In contrast to the eigenvalue approach, the data-driven approach to the estimation of the stabilization parameter circumvents the above issues as the computational effort and the memory footprint are independent of the cut configuration and, consequently, $n_{\text{ai}}$. Furthermore, the cost of the calculation of the stabilization parameter on a cutcell corresponds to the inference of the neural network model, which is expected to be lower than the effort associated with the eigenvalue approach given that the network size is relatively small. Another advantage of the data-driven estimate is the virtually gratuitous and non-intrusive support for GPUs out of the box thanks to the wide adaption of GPUs by machine learning frameworks, including TensorFlow, whereas the implementation of the eigenvalue approach for such devices would be non-trivial. Therefore, in the following we compare the performance of the eigenvalue approach on the CPU with the data-driven approach in terms of computational cost on both CPU and GPU devices.

The performance of the data-driven approach using the Intel Xeon E5-2630 v3 CPU and two modern GPUs, namely the NVIDIA Tesla V100 and the NVIDIA Ampere A100 is measured in terms of the total runtime for the estimation of the stabilization parameter on the entire training dataset, which includes 477,603 data points and shown in Figure~\ref{fig:performance_training}. We note that the tests are carried out in parallel. As network inference is typically dependent on the batch size of the inference data, the runtime with different batch sizes is reported. The optimal batch size was found to be 4,096 on the Intel Xeon, 32,768 on the NVIDIA V100 and 131,072 on the NVIDIA A100. We note that larger batch sizes did not seem to carry meaningful performance benefits. The general trend indicates that, as expected, larger batch sizes are more efficient and should, therefore, be preferred. The GPUs are observed to be up to around 42 times faster than the CPU using the optimal batch size on each device.

In order to compare the performance of the eigenvalue approach with the data-driven approach, as TensorFlow does not provide a straightforward method for single-thread testing, the lowerbound runtime of the eigenvalue approach as a function of $n_{\text{ai}}$ on the Intel Xeon is estimated by multiplying the single-thread runtime in Figure~\ref{fig:eigval_runtime} by the number of data-points in the training dataset and dividing it by the number of cores (16), assuming ideal parallel speedup. It is then possible to compare the ratio between the lowerbound runtime of the eigenvalue approach directly with the runtime of the data-driven approach in Figure~\ref{fig:performance_training} on the CPU as shown in Figure~\ref{fig:runtime_ratio_cpu} and on the GPU as shown in Figure~\ref{fig:runtime_ratio_gpu}. The data-driven approach demonstrates a clear advantage in terms of runtime, especially when $n_{\text{ai}}$ is large, up to $6.3 \times 10^{4}$ on the CPU and $2.7 \times 10^{6}$ on the GPU.

We focus on the integrated performance of the proposed data-driven approach next using a numerical benchmark. The physical domain consists in the shape shown in Figure~\ref{fig:benchmark}~(a), where the steady-state temperature field due to the boundary conditions shown in the same figure is sought. We note that the curved surface of the physical domain should create a representative variety of cut configurations. Adaptive mesh refinement, as shown in Figures~\ref{fig:benchmark}~(b) and (c), is used towards the boundary of the physical domain, resulting in a mesh with $988,837$ cells, of which $197,796$ are cutcells. It is important to highlight the high ratio between the number of cutcells and the total number of cells in the mesh as a consequence of adaptive mesh refinement. The data-driven estimate is utilized for all cutcells where the intersection between the physical boundary and the cell is either a line or can be closely approximated using a line according to the criterion explained in Section~\ref{sec:cut_config}; otherwise, the eigenvalue estimate is utilized. It is found that after adaptive refinement, the eigenvalue estimate did not need to be used for any cutcells, which is consistent with our observations for other examples, showcasing that the chosen representation of the cut configuration does not, in practice, hinder its usage. In order to evaluate the performance of the data-driven estimate in comparison with the eigenvalue estimate, the total computational cost of estimating the stabilization parameter for all cutcells in the domain as well as the quality of the estimate are considered in the following.
The total computational cost of the data-driven estimate for the given benchmark is measured at $2.84$ seconds on the Intel Xeon CPU. Given that the equivalent cost of the eigenvalue approach using 10 and 20 levels of integration levels is $8900$ and $28.74 \times 10^{5}$ seconds, respectively, see Figure~\ref{fig:eigval_runtime}, and given that the data-driven estimate is trained to match the quality of $n_{\text{ai}} = 20$, the data-driven approach offers substantial savings in terms of the total computational cost. The margin of such savings is even larger when the estimation of the stabilization parameter is offloaded to accelerators such as GPUs, as demonstrated in Figures~\ref{fig:performance_training} and~\ref{fig:runtime_ratio_gpu}. Furthermore, it is found that the difference between the stabilization parameter from the eigenvalue estimate and the data-driven estimate is capped at $5\%$, indicating that the quality of the estimation remains virtually the same. We note that the solution of finite cell problem remains identical down to machine accuracy between the eigenvalue and data-driven estimates.

\section{Conclusions}
\label{sec:conclusions}
In this work, we propose a data-driven approach to the estimation of the stabilization parameter in Nitsche's method. We use the finite cell formulation of Poisson's equation as model problem and compare the proposed data-driven estimate with the conventional eigenvalue-based estimate in terms of accuracy and computational efficiency. While the computational effort and memory footprint associated with the eigenvalue approach depends on the cut configuration, the data-driven estimate offers a more efficient alternative with $O(1)$ computational complexity. It is seen that the data-driven approach can accurately estimate the stabilization parameter for a wide range of cut configurations using a relatively small neural network. As a result, the data-driven estimate is by far more computationally efficient than its eigenvalue-based counterpart. The data-driven estimate can replace the eigenvalue estimate in the simulation pipeline in existing codes without major modifications. Furthermore, thanks to the wide adoption of accelerators by machine learning platforms, the data-driven estimate can be used on devices such as GPUs virtually out of the box with minimal implementation effort. The results from the numerical benchmarks show that the data-driven estimate can be many times more efficient compared to the eigenvalue estimate, while its relative error remains under $5\%$. The efficiency margin of the data-driven estimate is shown to be even larger using modern GPU devices such as the NVIDIA Tesla V100 and the NVIDIA Ampere A100.

\begin{acknowledgements}
Financial support was provided by the German Research Foundation ({\it Deutsche Forschungsgemeinschaft, DFG}) in the framework of subproject C4 of the Collaborative Research Center SFB 837 {\it Interaction Modeling in Mechanized Tunneling}, grant number 77309832. This support is gratefully acknowledged. We also gratefully acknowledge the computing time on the computing cluster of the SFB 837 and the Department of Civil and Environmental Engineering at Ruhr University Bochum, which was used for the presented numerical studies.
\end{acknowledgements}
\printbibliography

\end{document}